\documentclass[12pt]{amsart}

\newcommand\real{{\mathrm I}\!{\mathrm R} }

\newcommand\rat{{\mathrm Q}\kern-.65em {}^{{}_/ }}

\newtheorem{corollary}{Corollary}
\newtheorem{remark}{Remark}

\newtheorem{theorem}{Theorem}
\newtheorem{lemma}{Lemma}

\begin{document}
\title{Measures maximizing topological pressure}
\author{Abdelhamid Amroun}
\maketitle
\address{\begin{center}
Universit\'e Paris-Sud, D\'epartement de Math\'ematiques,
CNRS UMR 8628, 91405 Orsay Cedex France
\end{center}}
\address{}
\begin{abstract} We give a general method on the way of approximating
equilibrium states for a dynamical system of a compact metric
space. 
\end{abstract}
\section{Introduction}
We show in this paper how to buit equilibrium states for a dynamical
system on a compact metric space $X$ starting from the definition of the
topological pressure. For each finite set $E$ of ``sufficiently
separated'' points, we consider the probability measure $\mu_{E}$ with
support $E$. We then prove that the weak limits of $\mu_{E}$, as
$|E|\rightarrow \infty$, are equilibrium states. 
In fact, $|E|$ increases exponentially fast and $|E|\approx
e^{th_{top}}$, where $h_{top}$ is the topological entropy of the
dynamical system.  The result works for
flows as well as for maps of $X$. We begin with a general
and abstract result which gives a sufficient condition under which
the weak limits of a sequence of probability measures are contained in
a given closed and convex subset of $\mathcal{P}(X)$, the space of probability
measures of $X$. We apply this result to convex and lower semicontinuous
functionals on $\mathcal{P}(X)$. It results a way of approximating
zeros of such functionals. In the setting of dynamical systems, these
functionals are given by the Legendre transform of the topological
pressure considered as a functional on continuous potentials on $X$.
In certain cases, the result follows from a direct application of the
ergodic theorem. We give an example of such situation, by considering
probability measures supported on finite trajectories of a diffusion process
associated with a second order differential operators.   
\section{An abstract model and main result}
Let $X$ be a compact metric space and denote by $\mathcal{P}(X)$
the space of probability masures on $X$ equipped with the topology of
weak convergence of measures. 
Let $\beta_{t}$, $t\in \real$, be a probability measure on $X$ and
consider a random variable $\delta_{t}$ on the measure space
$(X,\beta_{t} )$ with values in $\mathcal{P}(X)$. Thus, for each $x\in
X$, $\delta_{t}(x)$ is a probability measure on $X$.

Consider the law $\nu_{t}$ of $\delta_{t}$ with 
respect to $\beta_{t}$, that is, $\nu_{t}:=\beta_{t} \circ \delta_{t}^{-1}$.
Also, the expectation of $\delta_{t}$ is given by the probability
measure 
\begin{equation}
\mu_{t}:=E_{\beta_{t}}(\delta_{t})=\int_{X}\delta_{t}(x)d\beta_{t}.
\end{equation}
We state the main result of this section.
\begin{theorem} Let $(\beta_{t})_{t\geq 0}$ be a family of probability
measures on $X$, and for each $t\geq 0$ a map $\delta_{t}$ on $(X,\beta_{t})$
with values in the space $\mathcal{P}(X)$. Suppose 
there is a nonempty closed and convex subset $\Gamma$ of $\mathcal{P}(X)$ such
that: for any open neighborhood $U\subset \mathcal{P}(X)$ of
$\Gamma$ we have, 
\begin{equation}
\lim_{t\rightarrow +\infty} \nu_{t}(U)=1.
\end{equation}
Then, any weak limit $\mu$ of
$\mu_{t}:=E_{\beta_{t}}(\delta_{t})$, as $t\rightarrow +\infty$ is
contained in $\Gamma$. In particular, if $\Gamma=\{m\}$ then $\mu_{t}$
converges to $m$. 
\end{theorem}
Note that the condition $(2)$ in the theorem is equivalent to assume
that $\lim_{t\rightarrow +\infty} \nu_{t}(U^{c})=0$ since
$\nu_{t}(U)+\nu_{t}(U^{c})=1$. This condition means
that the proportion of probability measures of the type $\delta_{t}$
which are close to $\Gamma$ is assymptotically $1$.
On the other hand, in the case where $\beta_{t}$ do not depend on $t$,
there are non trivial examples where the convegence of
$E_{\beta}(\delta_{t})$ to some measure is a 
consequence of the ergodic theorem (see section $4$).

The applications given below are based on the following special case
of the main Theorem $1$. Let $J: \mathcal{P}(X) \rightarrow [0,
\infty]$ be a lower semicontinuous functional and set $\Gamma :=\{J=0\}$.
This is a compact subset of
$\mathcal{P}(X)$ and furthermore, the convexity of $\Gamma$ follows from
the convexity of the functionl $J$.  

In practice, to  obtain $(2)$ one can prove a stronger result, namely a
large deviation upper bound for the process $\nu_{t}$: for any closed
subset $K$ of $\mathcal{P}(X)$,
\begin{equation}
\limsup_{t\rightarrow +\infty}\frac{1}{t}\log \nu_{t}(K) \leq
-\inf_{m\in K}J(m).
\end{equation}
A sufficient condition which imply this upper boud is given in
\cite{ki} (condition $(1.3)$ p506). In the notations of Theorem $1$,
this condition is equivalent to say that the limit
\begin{equation}
\lim_{t\rightarrow +\infty}
\frac{1}{t}\log E_{\beta_{t}}(e^{t\int_{X} fd\delta_{t}})
\end{equation}
exists for any continuous function $f$ on $X$. When this limit exists,
it coincides with the topological pressure $P(f)$ introduced in the
sequel. Furthermore, the functional $J$ in this case will be the
convex conjugate of $P$ in a certain sens (see section $3$).

Let then $U$ be an open neighborhood of $\Gamma :=\{J=0\}$.
Since $J$ is lower semicontinuous and the set $U^{c}$ is compact,
there will exist $m_{0}\in U^{c}$ such that $\inf_{m\in
  U^{c}}J(m)=J(m_{0})>0$. Thus, for $t$ sufficiently large we will
have $\nu_{t}(U) \geq 1-e^{-tJ(m_{0})}$, which implies $(2)$. We
summarize this discussion in the following result which is a corollary
of Theorem $1$.
\begin{theorem}
Let $(\beta_{t})_{t\geq 0}$ be a family of probability
measures on $X$, and for each $t\geq 0$ a map $\delta_{t}$ on $(X,\beta_{t})$
with values in the space $\mathcal{P}(X)$.
Suppose there exists a lower semicontinuous convex
functional $J: \mathcal{P}(X) \rightarrow [0, \infty]$ such that, for
any open neighborhood $U$ of $\{J=0\}$ we have,
\begin{equation}
\limsup_{t\rightarrow +\infty}\frac{1}{t}\log \nu_{t}(U^{c}) \leq
-\inf_{m\in U^{c}}J(m).
\end{equation}
Then, $J(\mu)=0$ for any weak limit $\mu$ of
$\mu_{t}:=E_{\beta_{t}}(\delta_{t})$, as $t\rightarrow +\infty$.
\end{theorem}
\subsection{Proof of the main result}
\begin{proof}
 We endow $\mathcal{P}(X)$ with a
distance $d$ compatible with this topology: 
take a countable base $\{g_{1}, g_{2}, \cdots\}$ of the separable space
$C_{\real}(X)$, where $\|g_{k}\|=1$ for all $k$, and set:
\[
d(m, m'):=\sum_{k=1}^{\infty}2^{-k}\left |\int g_{k}dm - \int g_{k}dm'
\right |.
\] 
Let $V\subset \mathcal{P}(X)$ be a convex open 
neighborhood of $\Gamma$ and $\gamma>0$.
Consider a finite open cover $(B_{i}(\gamma))_{i\leq N}$ of
$\Gamma$ by balls of diameter $\gamma$ all
contained in $V$ (if $B_{i}(\gamma)$ is not entirely contained in $V$
we just take its restriction to $V$).  

Decompose the set $U:=\cup_{i=1}^{N}B_{i}(\gamma)$ as follows,
\[
U=\cup_{j=1}^{N'}U_{j}^{\gamma},\ for \ some \ N'\geq N,
\]
where the sets $U_{j}^{\gamma}$ (which are not necessarily open) are
disjoints and contained in one of the balls $(B_{i}(\gamma))_{i\leq N}$.
We have 
\[
\Gamma \subset U \subset V,
\] 
and $\sum_{j=1}^{N'}\nu_{n}(U_{j}^{\gamma})=\nu_{n}(U)$.
We fix in each $U_{j}^{\gamma}$ a probability measure $p_{j}$,
$j\leq N'$, and let $p_{0}$ be a probability  measure distinct from the above
ones (for example take $p_{0} \in V\backslash U$). Set,
\[
\omega_{t}:=\sum_{j=1}^{N'}p_{j}1_{\delta_{t}^{-1}(U_{j}^{\gamma})}+
(1-\nu_{t}(U))p_{0}.
\]
We have,
\begin{equation}
E_{\beta_{t}}(\omega_{t})=
\sum_{j=1}^{N'}\nu_{t}(U_{j}^{\gamma})p_{j}+
(1-\nu_{t}(U))p_{0}. 
\end{equation}
The probability measure $E_{\beta_{t}}(\omega_{t})$ are contained in $V$ since
it is a convex combination of elements of the convex set $V$. 
We have then
$d(\mu_{t}, V)\leq d(\mu_{t},E_{\beta_{t}}(\omega_{t}))$. 
We will show that
\[
d(\mu_{t},E_{\beta_{t}}(\omega_{t}))\leq \gamma
\nu_{t}(U)+ 3\nu_{t}(U^{c}), 
\]
where $U^{c}=\mathcal{P}(X)\backslash U$.

Consider the finite positive measures,
\[
\mu_{t,V}:=E_{\beta_{t}}\left ((1_{V}\circ
\delta_{t})\delta_{t} \right ).
\]Evaluated on a function $g$ this gives,
\[
\mu_{t,V}(g):=E_{\beta_{t}}\left ((1_{V}\circ
\delta_{t})\delta_{t}(g) \right )=
\int_{X}1_{V}(\delta_{t}(x))
\left ( \int_{X}g(y)d\delta_{t}(x)(y)\right )d\beta_{t}(x).
\]
Note that for any continuous function $g$ on $X$ with $\|g\|=1$ and
any subset $E$ of $\mathcal{P}(X)$ we
have 
\begin{equation}
|\mu_{t,E}(g)|\leq E_{\beta_{t}}(1_{E}\circ
\delta_{t})=\nu_{t}(E).
\end{equation}
By the definition of $d$ we have to estimate:
\[
d(\mu_{t},E_{\beta_{t}}(\omega_{t}))=
\sum_{k=1}^{\infty}2^{-k}\left |\mu_{t}(g_{k})- 
E_{\beta_{t}}(\omega_{t})(g_{k}) \right |.
\]
First write,
\begin{eqnarray*}
&& d(\mu_{t},E_{\beta_{t}}(\omega_{t}))\\
&\leq&\sum_{k\geq 1}2^{-k}|\mu_{t}(g_{k}) - \mu_{t,V}(g_{k})|+
\sum_{k\geq 1}2^{-k}|\mu_{t,V}(g_{k}) -
E_{\beta_{t}}(\omega_{t})(g_{k})|.
\end{eqnarray*}
From $(6)-(7)$, the definition of $\mu_{t}$ and $\mu_{t,V}$ and
the fact that $U\subset V$, we get
\begin{equation}
\sum_{k\geq 1}2^{-k}|\mu_{t}(g_{k}) - \mu_{t,V}(g_{k})| 
=\sum_{k\geq 1}2^{-k}|\mu_{t,V^{c}}(g_{k})| 
\leq \nu_{t}(U^{c}). 
\end{equation}
It remains to show that 
\[
\sum_{k\geq 1}2^{-k}|\mu_{t,V}(g_{k}) -
E_{\beta_{t}}(\omega_{t})(g_{k})|\leq \gamma
\nu_{t}(U)+\nu_{t}(U^{c}).
\] 
We have
\[
\sum_{k\geq 1}2^{-k}\left
|\mu_{t,V}(g_{k})-E_{\beta_{t}}(\omega_{t})(g_{k})\right |=
\sum_{k\geq 1}2^{-k}\left|A-B\right |
\]
where,
\begin{eqnarray*}
A&:=&\sum_{j=1}^{N'}\left (E_{\beta_{t}}\left ((1_{U_{j}^{\gamma}}\circ
\delta_{t})\delta_{t}\right )(g_{k})
-\nu_{t}(U_{j}^{\gamma})p_{j}(g_{k})
\right ),\ and \\
B&:=&E_{\beta_{t}}\left ((1_{V\backslash U}\circ
\delta_{t})\delta_{t}\right )(g_{k})
+(1-\nu_{t}(U))p_{0}(g_{k}).
\end{eqnarray*}
We have for all $k\geq 1$,
\begin{equation}
\sum_{k\geq 1}2^{-k}A\leq \gamma
\sum_{j=1}^{N'}\nu_{t}(U_{j}^{\gamma})=\gamma \nu_{t}(U). 
\end{equation}
and
\begin{equation}
\sum_{k\geq 1}2^{-k}B\leq \nu_{t}(V\backslash
U)+(1-\nu_{t}(U)) 
\leq 2\nu_{t}(U^{c}).
\end{equation}
Putting toguether $(8)-(10)$ we get,
\[
d(\mu_{t},E_{\beta_{t}}(\omega_{t}))
\leq \gamma \nu_{t}(U)+3\nu_{t}(U^{c}). 
\]
This implies the desired inequality,
\[
d(\mu_{t},V)\leq \gamma
\nu_{t}(U)+3\nu_{t}(U^{c}). 
\]
Now, since $U$ is open, we have by assumption, $\lim_{t\rightarrow
\infty}\nu_{t}(U)=1$ and $\lim_{t\rightarrow
\infty}\nu_{t}(U^{c})=0$. Thus, $\limsup_{t\rightarrow
\infty}d(\mu_{t},V)\leq \gamma$, for all $\gamma>0$. We
conclude that $\limsup_{t\rightarrow
\infty}d(\mu_{t},V)=0$. The neighborhood $V$ of
$\Gamma$ being arbitrary, this implies that the all limit
measures of $\mu_{t}$ are contained in $\Gamma$. 
\end{proof}
\begin{remark}
It is not necessary for $V$ to be convex. Indeed,
let $V$ be an open neighborhood of $\Gamma$, and consider
finite open cover $(B_{i}(\gamma))_{i\leq N}$ of
$\Gamma$ by balls of diameter $\gamma$ all contained in
$V$ where each ball $B_{i}(\gamma)$ is centred at a measure
$m_{i} \in \Gamma$. Fix also a measures $m_{0} \in
\Gamma$ distinct from the above ones (since
$\Gamma$ is convex, if it contains two elements it
contains the segment formed by them). Now define
$E_{\beta_{t}}(\delta_{t})$ 
as in $(6)$. 
The difference here resides in the fact that each $m_{i}$ can
appear more than one time in $(6)$.
However, $E_{\beta_{t}}(\delta_{t})$ is still a
convex combination of elements in $\Gamma$, so that
$E_{\beta_{t}}(\delta_{t}) \in \Gamma$. 
\end{remark}

\section{Application to dynamical systems}
Let us explain how to apply the above results in the context of dynamical
systems. Leaving the details and the precise statements for later, we
begin by a general description of the method. 
Let $X$ be a compact metric space and $\phi: \real \times X
\rightarrow X$ a
continuous (smooth) flow acting on it.
The same can be done for maps of $X$. 
Set $\phi^{t}:=\phi (\cdot, t):
X\rightarrow X$.
The functional $J$ has then the following form,
\begin{equation}
J(\mu)=\sup_{\omega \in C_{\real}(X)} \left (\int \omega d\mu-
Q(\omega)\right ),
\end{equation}
where $Q: C_{\real}(X) \rightarrow \real$ is continuous and convex,
$\phi$-invariant i.e, 
$Q(\omega \circ \phi^{t})=Q(\omega)$ for $t\in \real$, and $Q(0)=0$. 
Thus $J\geq 0$ and $J$ is lower semicontinuous. Since for any
$\omega$, the map $\mu \rightarrow (\int \omega d\mu -Q(\omega))$ 
is convex, we deduce by duality  that,
\begin{equation}
Q(\omega)=\sup_{\mu \in \mathcal{P}(X)} \left (\int \omega d\mu-
J(\mu)\right ).
\end{equation}
We are interested in the subset $\mathcal{P}_{0}(X)$ of
$\mathcal{P}(X)$ of 
probability measures $\mu$ such that $J(\mu)=0$,
\begin{equation}
\mathcal{P}_{0}(X):=\{\mu \in \mathcal{P}(X): J(\mu)=0\}. 
\end{equation}
Elements of $\mathcal{P}_{0}(X)$ which are invariant by $A$, are called
equilibrium states of the dynamical system. 
Theorem $2$ in this case says that, if we have
the upper bound $(5)$, then $\mu \in \mathcal{P}_{0}(X)$ for any weak
limit $\mu$ of $\mu_{t}=E_{\beta_{t}}(\delta_{t})$. 

The objective in what follows is to specify all the ingredients of
Theorem $2$ and to establish $(5)$. 

A probability mesure $\mu$ on $X$ is said to be invariant by
$\phi$, if for all continuous function $\omega$ and all $t\in \real$
we have $\int_{X} \omega (\phi^{t}(x)) 
d\mu=\int_{X} \omega (x) d\mu$. 
Let $\mathcal{P}_{inv}(X)$ be the closed and convex subspace of
$\mathcal{P}(X)$ which are invariant by the flow $\phi$. 

\subsection{Topological pressure and the functional $J$}
The topological pressure $P(f)$ of a continuous function
$f :X \rightarrow X$ is defined using separated sets of
$X$ (see \cite{walt}) as follows:
\begin{equation}
P(f)=\lim_{\epsilon \rightarrow 0}\limsup_{t\rightarrow +\infty}
\frac{1}{t}\log\sup_{E}\sum_{x\in E}e^{\int_{0}^{t}f(\phi^{s}x)ds},
\end{equation}
where the $\sup$ is taken over the $(t, \epsilon)$-separated sets
$E\subset X$. Since $X$ is a compact, in formula $(14)$ we can take the
$\sup$  over maximal separated sets \cite{walt}. There is also an equivalent
definition based on spanning sets \cite{walt}.

It is well known that $P(f)$ satisfies the following variational
principle \cite{walt}, 
\begin{equation}
P(f)=\sup_{\mu \in \mathcal{P}_{inv}(X)}\left (h(\mu)+\int fd\mu \right ).
\end{equation}
For $f=0$ formula $(12)$ gives the usual variational principle for the
topological entropy $h_{top}$ of the flow 
$\phi$,
\begin{equation}
h_{top}=\sup_{\mu \in \mathcal{P}_{inv}(X)} h(\mu).
\end{equation}
Let $\mathcal{P}_{e}(f)$ be the subset of measures $\mu \in
\mathcal{P}_{inv}(X)$ realizing the $\sup$ in $(15)$. These maximazing
measures are called equilibrium states corresponding to the potential $f$.

We assume in the sequel that the entropy map $m\rightarrow h(m)$ is
upper semicontinuous (u.s.c). Then $h_{top} <+\infty$ and
$\mathcal{P}_{e}(f)$ is a nonempty closed and convex
subset of $\mathcal{P}_{inv}(X)$ (see \cite{walt}). For example, by an
important result of Newhous \cite{new}, if $X$ is a compact manifold
equipped with a $C^{\infty}$ Riemannian metric then the entropy map is u.s.c.
In this context, the main example of flow we have in mind is the
geodesic flow \cite{am}.  

Define the functional $Q_{f}$ by,
\begin{equation}
Q_{f}(\omega):=P(f+\omega)-P(f),\ for \ \omega \in C_{X}(\real).
\end{equation}
The functional $J_{f}$ is then defined on the space of
probability measures on $X$ by,
\begin{equation}
J_{f}(\mu):=\sup_{\omega}\left (\int \omega d\mu-Q_{f}(\omega) \right ),
\end{equation}
where the $\sup$ is over the continuous finctions $\omega$ on
$X$. Since $Q_{f}(0)=0$, we have $J_{f}\geq 0$.
In \cite{am} (Lemma 1,2) we proved,
\begin{lemma} 
\noindent
\begin{enumerate}
\item 
$Q_{f}$ is $\phi$-invariant, that is
$Q_{f}(\omega \circ \phi_{t})=Q_{f}(\omega)$ for all
continuous function $\omega$ and $t\in \real$. Moreover,
$Q_{f}$ is convex and continuous on continuous functions.  
$J_{f}$ is convex and lower semicontinuous functional and by duality,
\[
Q_{f}(\omega)=\sup_{\mu \in \mathcal{P}(X)}\left (\int \omega
d\mu-J_{f}(\mu) \right ).
\]
\item Set for any invariant probability measure $\mu$,
\[
I(\mu):=P(f)-\left (h(\mu)+\int_{X}f d\mu \right ).
\]
Then $Q_{f}(\omega)=\sup_{\mu \in \mathcal{P}_{inv}(X)}\left (\int \omega
d\mu-I(\mu)\right )$. In other words, the functionals $I$ and $J_{f}$
agree on invariant measures.
\item We have,
\[
\mathcal{P}_{e}(f)=\{J_{f}=0\}\cap \mathcal{P}_{inv}(X).
\]
\item The uniqueness of the solution of the equation $J_{f}=0$ in
 $\mathcal{P}(X)$ is equivalent the uniqueness
  of the equilibrium state corresponding to $f$.
\end{enumerate}
\end{lemma}
Note that a probability measure $m$ for which we have $J_{f}(m)=0$ is
not necessarily invariant, however the set $\{J_{f}=0\}$ is invariant
by the flow.
Given $x\in X$ and $t>0$, we define a probability measure
$\delta_{t}(x)$ on $X$ by:
\begin{equation}
\int \omega d\delta_{t}(x):=\frac{1}{t}\int_{0}^{t}\omega(\phi^{s}x)ds.
\end{equation}
The main result of this section is,
\begin{theorem} For any real continuous function $f$ on $X$ we have:
\begin{enumerate}
\item Let $\gamma >0$. There exists $\epsilon_{0}>0$ such that for all
$\epsilon <\epsilon_{0}$, there exists a sequence $(E_{n},
  t_{n})_{n\geq 0}$ where $E_{n}$ is a $(t_{n}, \epsilon)$-separated
  set in $X$ and $ t_{n}\geq n$ such that,
\[
\sum_{x\in E_{n}}e^{\int_{0}^{t_{n}}f(\phi^{s}x)ds}\geq
e^{t_{n}(P(f)-\gamma)}, \ for \ all \ n\geq 0.
\]
\item Let $\mu_{n}$ be the probability measures defined by, 
\[
\mu_{n}:=\mu_{t_{n}}=
\frac{\sum_{x\in E_{n}}e^{\int_{0}^{t_{n}}f(\phi^{s}x)ds}\delta_{t_{n}}(x)}
{\sum_{x\in E_{n}}e^{\int_{0}^{t_{n}}f(\phi^{s}x)ds}},
\]
where $(t_{n},E_{n})$ is the sequence determined in $(1)$.
Then any weak limit $\mu$ of $\mu_{n}$, as $n\rightarrow \infty$ is
invariant by the flow and satisfies $J_{f}(\mu)=0$. The measure $\mu$
is then an equilibrium state. 
\end{enumerate}
\end{theorem}
As the proof will show, the conclusion in part $(2)$ of Theorem $3$ is true
for any sequence $E_{n}$ of $(t_{n}, \epsilon)$-separated sets in $X$
for which Theorem 3 (1) holds. We describe now examples of such situations.

Consider a smooth compact Riemannian manifold $M$ with negative
curvature. The set of primitive closed geodesic (which represent
different free homotpy class) of period $\leq t$ is $(\epsilon,
t)$-separated whenever $\epsilon < inj(M)$, 
where $inj(M)$ is the injectivity radius of $M$. Moreover,
condition $(1)$ of Theorem $3$ holds here since the exponential growth
rate of these geodesics is given by the topological pressure \cite{bowen}
\cite{bowe} \cite{boru} \cite{led},
\[
P(f)=\lim_{t\rightarrow +\infty}\frac{1}{t}\log
\sum_{c:\ l(c)\leq t}e^{\int_{c}f},
\]where the sum is over closed geodesics (primitive) of length at
most $t$. The precise statement is as follows.
\begin{corollary}
Let $M$ be a compact manifold equipped with a $C^{\infty}$ Riemannian
metric of negative curvature. Let $\phi^{t}: X\rightarrow X$ be the
geodesic flow acting on the unit tangent bundle $X=T^{1}M$ of
$M$. Then, for any H\"older continuous function $f: X\rightarrow
\real$, the probability measures defined by
\[
\mu_{t}(\omega):=\frac{\sum_{c:\ l(c)\leq t}e^{\int_{c}f}
\left ( \frac{1}{l(c)}\int_{c}\omega \right )}
{\sum_{c:\ l(c)\leq t}e^{\int_{c}f}},
\]
converge to the unique equilibrium state of the geodesic flow
corresponding to $f$.
 \end{corollary}

\subsection{Proof of Theorem $3$}
\subsubsection{Proof of Part $(1)$}
\begin{proof}
The proof is a consequence of the definition $(1)$. Indeed, let
$\eta>0$ be fixed. There exists $\epsilon_{0}>0$ such that for all
$0<\epsilon <\epsilon_{0}$:
\[
\inf_{n\geq 0}\sup_{t\geq n}\frac{1}{t}\log 
\sup_{E:(t,\epsilon)-separated}\sum_{x\in
  E}e^{\int_{0}^{t}f(\phi^{s}x)ds}\geq P(f)-\eta.
\]
From this we deduce that for each $n\geq 0$ there exists $t_{n}\geq n$
and an $(t_{n},\epsilon)$-separated set $E_{n}$ such that,
\[
\sum_{x\in
  E_{n}}e^{\int_{0}^{t_{n}}f(\phi^{s}x)ds}\geq
  e^{t_{n}(P(f)-\eta)}. 
\]
\end{proof}
\subsubsection{Proof of Part $(2)$}
\begin{proof}
Define the atomic probability measures with support in $E_{n}$:
\begin{equation}
\beta_{n}(B)=\beta_{t_{n}}(B):=\frac{\sum_{x\in B\cap
E_{n}}e^{\int_{0}^{t_{n}}f(\phi^{s}x)ds}}{\sum_{x\in 
E_{n}}e^{\int_{0}^{t_{n}}f(\phi^{s}x)ds}}.
\end{equation}
Then, the measure $\delta_{t_{n}}$
defined in $(19)$, can be seen as random variable on the measure space
$(X,\beta_{n})$ and the measures $\mu_{n}$ in Theorem $3$ as the
  expectation $\mu_{n}=E_{\beta_{n}}(\delta_{t_{n}})$. 
Consider the image $\nu_{n}$ of $\beta_{n}$ under the map
$\delta_{t_{n}}$:
\begin{equation}
\nu_{n}:=\beta_{n} \circ \delta_{t_{n}}^{-1}.
\end{equation}
To prove part $(2)$ we have to check the upper bound in Theorem $2$
for $\{\nu_{n}\}$ namely, for any closed subset $K$ of
$\mathcal{P}(X)$, 
\[
\limsup_{n\rightarrow +\infty}\frac{1}{t_{n}}\log \nu_{n}(K)\leq
-J_{f}(K):=-\inf_{m\in K}J_{f}(m). 
\]
We follow \cite{pol}, \cite{ki} and \cite{am}.
Let $\eta>0$. Observe that the compact set $K$ is contained the union of
open sets,
\[
K\subset \cup_{\omega}\{\mu \in \mathcal{P}(X) :\int
\omega dm-Q_{f}(\omega)>J_{f}(K)-\eta\}.
\]
There exists then a finite number of continuous
functions $\omega_{1}, \cdots, \omega_{l}$ such that $K\subset
\cup_{i=1}^{l}K_{i}$, where 
\[
K_{i}=\{m\in \mathcal{P}(X): \int
\omega_{i}dm-Q_{f}(\omega_{i})>J_{f}(K)-\eta\}. 
\]
We have $\nu_{n}(K)\leq \sum_{i=0}^{l}\nu_{n}(K_{i})$ where,
\begin{equation}
\nu_{n}(K_{i})=
\frac{\sum_{x\in E_{n}:\delta_{t_{n}}(x)\in K_{i}}
e^{\int_{0}^{t_{n}}f(\phi^{s}x)ds}}{\sum_{x\in 
E_{n}}e^{\int_{0}^{t_{n}}f(\phi^{s}x)ds}}.
\end{equation}
For $n$ sufficiently large we have,
\begin{eqnarray*}
&&\sum_{x\in E_{n}:\delta_{t_{n}}(x)\in K_{i}}
e^{\int_{0}^{t_{n}}f(\phi^{s}x)ds}\\
&\leq& \sum_{x\in E_{n}:\delta_{t_{n}}(x)\in K_{i}}
e^{\int_{0}^{t_{n}}f(\phi^{s}x)ds}
e^{t_{n}(\int
 \omega_{i}d\delta_{t_{n}}(x)-Q_{f}(\omega_{i})-(J_{f}(K)-\eta))}\\   
&\leq&e^{t_{n}(-Q_{f}(\omega_{i})-(J_{f}(K)-\eta))}
\sum_{x\in E_{n}}
e^{\int_{0}^{t_{n}}(f+\omega_{i})(\phi^{s}x)ds}.
\end{eqnarray*}
By the definition of the pressure and the fact that $t_{n}\geq
n$ we also have, 
\begin{eqnarray*}
&&\limsup_{n\rightarrow \infty}\frac{1}{t_{n}}\log 
\sum_{x\in E_{n}}e^{\int_{0}^{t_{n}}(f+\omega_{i})(\phi^{s}x)ds}\\
&\leq&\limsup_{n\rightarrow \infty}\frac{1}{t_{n}}\log 
\sup_{E:(t_{n},\epsilon)-separated}\sum_{x\in
  E}e^{\int_{0}^{t_{n}}(f+\omega_{i})(\phi^{s}x)ds}\\ 
&\leq& \limsup_{t\rightarrow \infty}\frac{1}{t}\log 
\sup_{E:(t,\epsilon)-separated}\sum_{x\in
  E}e^{\int_{0}^{t}(f+\omega_{i})(\phi^{s}x)ds}\\
&=&P(f+\omega_{i}).
\end{eqnarray*}
Thus, for $n$ sufficiently large we obtain,
\[
\sum_{x\in E_{n}}e^{\int_{0}^{t_{n}}(f+\omega_{i})(\phi^{s}x)ds}
\leq e^{t_{n}(P(f+\omega_{i})+\eta)}.
\]
Consequently we obtain for $\nu_{n}(K)$,
\begin{eqnarray*}
\nu_{n}(K)&\leq& 
\sum_{i=1}^{l}\nu_{n}(K_{i})\\ 
&\leq&
\sum_{i=1}^{l}\left
(e^{t_{n}(-P(f)+\eta-Q_{f}(\omega_{i})-(J_{f}(K)-\epsilon))}  
\sum_{x\in
  E_{n}}e^{\int_{0}^{t_{n}}(f+\omega_{i})(\phi^{s}x)ds}\right)\\
&\leq& \sum_{i=1}^{l}\left
(e^{t_{n}(-P(f)+\eta-Q_{f}(\omega_{i})-(J_{f}(K)-\eta))}
 e^{t_{n}(P(f+\omega_{i})+\eta)} \right).
\end{eqnarray*}
But by $(17)$ we have
$P(f+\omega_{i})-P(f)=Q_{f}(\omega_{i})$. We deduce finally
that,
\[
\nu_{n}(K)\leq e^{t_{n}(-J_{f}(K)+2\eta))}
\]
Take the logarithme, divide by $t_{n}$ and the $\limsup$,
\[
\limsup_{n\rightarrow \infty}\frac{1}{t_{n}}\log \nu_{n}(K) 
\leq -J_{f}(K)+2\eta.
\]
$\eta$ being arbitrary, this proves part $(2)$.
\end{proof}
\section{Second order differential operators}
The main futers of this section is to give an example of what happend
if the measures $\beta_{t}$ do not depend on $t$ (Theorem 1). As we
will see, the convergence of $\mu_{t}$ towards the unique equilibrium
state is a consequence of the ergodic theorem. 
The materials of this section are taken from \cite{ki} and
\cite{ki1}. 

Let $M$ be a locally compact finite dimensional manifold
$M$ equipped with a $C^{2}$ Riemannian metric and $G$ a connected open
subset of $M$ with smooth boundary $\partial G$ such that
$\overline{G}=G\cup \partial G$. 
Let $L$ be a second order elliptic differential operator with
$C^{2}$ coefficients on $M$ if $G=M$ (i.e $M$ compact) or in some
neighborhood of $\overline{G}$ if $G\ne M$ and $L1=0$. This last
condition means that $L$ has locally the form
\[
L=\sum_{ij} a^{ij}(x)\frac{\partial^{2}}{\partial x_{i}\partial x_{j}}
+\sum_{i} b^{i}(x)\frac{\partial}{\partial x_{i}}.
\]
We denote by $\mathcal{P}(\overline{G})$ the space of probability
measures on $\overline{G}$. 
For any continuous function $V$ on $M$ define the operator
\[
L_{V}:=L+V.
\]
Let $\Sigma (V)$ be the spectrum of $L_{V}$ corresponding to the
Dirichlet boundary conditions on $\partial G$ (with no boundary
conditions if $G=M$ is compact). Consider,
\[
\lambda_{V}:=\sup \{Re (\lambda): \lambda \in \Sigma (V)\}. 
\]
Then
$\lambda_{V}$ is a spectral value for $L_{V}$ i.e $\lambda_{V} \in
 \Sigma (V)$ and by \cite{dv1, dv2} on hase the representation
\begin{equation}
\lambda_{V}=\sup_{\mu \in \mathcal{P}(\overline{G})}\left ( \int Vd\mu
-I(\mu)\right ) 
\end{equation}
where $I(\mu)$ is the entropy of the measure $\mu$ given by
\begin{equation}
I(\mu)=-\sup_{u \in D_{+}(L)}\int \frac{Lu}{u}d\mu
\end{equation}
and $D_{+}(L)$ is the set of functions from the domain
of $L$ having positive lower bounds.

For H\"older continous function $V$ it is well known \cite{kr}
that $\lambda_{V}$ is an eigenvalue of $L_{V}$, in fact it is the
leading eigenvalue. The corresponding eigenspace is one dimensional
and the corresponding eigenfunction $r_{V}$ is positive. Moreover
\cite{ki1}, there exists a unique 
probability measure $\mu_{V}$ on $\overline{G}$ solving the variational
principle $(23)$,
\begin{equation}
\lambda_{V}= \int Vd\mu_{V}-I(\mu_{V}).
\end{equation}
Let $X_{V}(t)$ be the difusion process with values in $M$ and 
generator $L_{V}$ \cite{friedman}, defined on a
probability space $(\Omega, \mathcal{F}, P)$.
Consider the semigroup of operators (Perron-Frobenius operators)
acting on $C(\overline{G})$ 
\begin{equation}
T_{V}(t)g(x):=E_{x}\left ( 1_{\tau_{G}>t}g(X_{V}(t))
e^{\int_{0}^{t}V(X_{V}(s))ds}\right )
\end{equation}
where $E_{x}$ denotes the expectation for $X_{V}(0)=x$ and
\[
\tau_{G}:=\inf\{t\geq 0: X_{V}(t) \notin G\}.
\]
Set $P_{x}(A)=E_{x}(1_{A})$.

If $G=M$ is compact just set $\tau_{G}=\infty$. In case of
H$\ddot{o}$lder continous function $V$ then $e^{\lambda_{V}}$ is
an eigenvalue of $T_{V}(1)$ and in fact its leading eigenvalue and we
have
\begin{equation}
\lambda_{V}=\lim_{t\rightarrow +\infty}\frac{1}{t}\log \left (
T_{V}(1) \right ),\ where \ 1(x)=1.
\end{equation}
Furthermore, the measure $\mu_{V}$ maximize $(20)$ if and only if it
is an invariant measure for the Markov process with transition
operators
\[
\mathcal{T}_{V}(t)g(x):=e^{-\lambda_{V}}r_{V}(x)^{-1}T_{V}(t)(gr_{V})(x).
\]
Then for any continuous function $u$ on $\overline{G}$
\[
\int_{\overline{G}}\mathcal{T}_{V}(t)u(x)d\mu_{V}(x)=
\int_{\overline{G}}u(x)d\mu_{V}(x).
\]
Thus for H$\ddot{o}$lder continous potentials $V$, the maximizing
measure $\mu_{V}$ is the the only invariant measure for the Markov
process above. $\mu_{V}$ is then ergodic.

Consider now the random probability measures supported on the finite
trajectories of the diffusion $X_{V}(t)$
\begin{equation}
\delta_{t}:=\frac{1}{t}\int_{0}^{t}\delta_{X_{V}(s)}ds
\end{equation}
where $\delta_{X_{V}(s)}$ is the Dirac measure at
$X_{V}(s)$. By the ergodic theorem, for  $\mu_{V}$ a.e $x\in M$, and
$P_{x}$ a.e $\omega \in \Omega$, the probability measures $\delta_{t}$
converge wealky towards $\mu_{V}$. Thus, it is easy to see that the measures
$E_{x}(\delta_{t})=\int \delta_{t} dP_{x}$ converge towards $\mu_{V}$.

Now it follows from \cite{ki} that the measures $\nu_{t}:=P_{x}\circ
\delta_{t}^{-1}$ satisfy a large deviation principle in particular the
following upper bound : for any $x\in G$ and any closed subset $K\in
\mathcal{P}(\overline{G})$ we have
\[
\limsup_{t\rightarrow +\infty}\frac{1}{t}\log
P_{x}(\delta_{t} \in K)\leq -\inf_{m\in K}J(\mu)
\]
where $J(m)=J_{V}(m):=I(m)-\int Vdm+ \lambda_{V}$. Note that for $V\equiv 0$
this is just the conclusion of Theorem $4.1$ from \cite{ki}. The
functional $J$ is non negative and $J(m)=0$ if and only if $m=\mu_{V}$ (recall
that $V$ is H$\ddot{o}$lder continous). Since $\mu_{V}$ is unique,
Theorem $2$ tell us that $E_{x}(\delta_{t})=\int \delta_{t} dP_{x}$
converge towards $\mu_{V}$, which is the statement proved above.

\end{document}